\documentclass[12pt,a4paper]{article}

\usepackage[dvips]{graphics}
\usepackage[dvips]{epsfig}

\newtheorem{lemma}{Lemma}

\newtheorem{prop}[lemma]{Proposition}

\newcommand{\dimo}{\noindent \emph{Proof. }}
\newcommand{\qed}{\\ \rightline{$\Box$ \ \ \ \ \ \ \ \ \ \ \ \ \ \ \ }\\}

\usepackage{latexsym}
\usepackage{amsfonts}
\usepackage{amssymb}
\usepackage{amsmath}
\usepackage[english]{babel}

\begin{document}
\title{
A CATALOGUE OF ORIENTABLE 3-MANIFOLDS \\
TRIANGULATED BY $30$ COLOURED TETRAHEDRA \footnote{Work performed
under the auspicies of the G.N.S.A.G.A. of the C.N.R. (National
Research Council of Italy) and financially supported by M.U.R. of
Italy (project ``Propriet\`a geometriche delle variet\`a reali e
complesse'').}}

\author{Maria Rita CASALI  - Paola CRISTOFORI \\
Dipartimento di Matematica Pura ed Applicata \\
       Universit\`a di Modena e Reggio Emilia \\
       Via Campi 213 B \\  I-41100 MODENA (Italy)}

\maketitle

\abstract {The present paper follows the computational approach to
3-manifold classification via edge-coloured graphs, already
performed in \cite{[L]} (with respect to orientable 3-manifolds up
to 28 coloured tetrahedra), in \cite{[C$_1999$]} (with respect to
non-orientable 3-manifolds up to 26 coloured tetrahedra), in
\cite{[C$_1989$]} and \cite{[BGR]} (with respect to genus two
3-manifolds up to 34 coloured tetrahedra): in fact, by  automatic
generation and analysis of suitable edge-coloured graphs, called
{\it crystallizations}, we obtain a catalogue of all orientable
3-manifolds admitting coloured triangulations with 30 tetrahedra.
These manifolds are unambiguously identified via JSJ
decompositions and fibering structures.

It is worth noting that, in the present work, a suitable use of
elementary combinatorial moves yields an automatic partition of
the elements of the generated crystallization catalogue into
equivalence classes, which are proved to be in one-to one
correspondence with the homeomorphism classes of the represented
manifolds. } \endabstract

  \bigskip
  \par \noindent
  {\bf Key words}: orientable 3-manifold, crystallization, coloured
  triangulation, complexity.
\smallskip
  \par \noindent
 \smallskip
  \par \noindent
  {\bf AMS Mathematics Subject Classifications (1991)}:
   57Q15 - 57M15 - 57N10.
  \bigskip

\section{\hskip -0.7cm . Introduction}

\bigskip

Within the study of the topology of PL-manifolds a great attention
has been recently reserved to combinatorial representation
methods, enabling to produce and study (possibly with the aid of
suitable computer programs) exhaustive catalogues of ``small"
manifolds, with respect to a given ``complexity" criterion: let us
recall, as a very significant example, successive results about
closed orientable irreducible 3-manifolds whose minimal special
spines have increasing number of vertices, up to $11$
(\cite{[M$_1$]}, \cite{[O]}, \cite{[MP$_1$]}, \cite{[M$_3$]}), or
the analogous studies about closed non-orientable $\mathbb
P^2$-irreducible 3-manifolds whose minimal special spines have at
most $10$ vertices (\cite{[AM_1]}, \cite{[C$_2004$]},
\cite{[B$_1$]}, \cite{[B$_2$]}, \cite{[AM_2]}, \cite{[B$_3$]},
\cite{[B$_4$]}).

During the last thirty years, another representation theory for
PL-manifolds has been developed. Its principal feature is
generality, i.e. it can represent the whole class of piecewise
linear (PL) manifolds, without assumptions about dimension,
connectedness, orientability, irreducibility, $\mathbb
P^2$-irreducibility or boundary properties: see \cite{[P]},
\cite{[FGG]}, \cite{[BM]}, \cite{[Co]}, \cite{[V]}, \cite{[KL]},
or \cite{[BCG]} for a survey on the so-called {\it crystallization
theory}, which makes use of edge-coloured graphs (named also {\it
crystallizations}, under suitable conditions) as a representation
tool.

Note that, in virtue of the purely combinatorial nature of the
representing objects, crystallization theory turns out to be
particularly suitable to computer enumeration.  From this
view-point, the main existing results concerning the whole class
of closed orientable 3-manifolds are described in Lins's book
\cite{[L]}, where a catalogue of all orientable 3-manifolds
represented by crystallizations with at most 28 vertices is
produced and analyzed.\footnote{See also \cite{[CC]}, where an
unambiguous identification of all elements of Lins's catalogue is
given, through JSJ decompositions and fibering structures.} On the
other hand, \cite{[C$_1999$]} takes into account the
non-orientable case, while other works deal with restricted
classes of 3-manifolds (for example, both orientable and
non-orientable euclidean 3-manifolds in \cite{[Va]}, genus two
orientable 3-manifolds in \cite{[C$_1989$]} and \cite{[BGR]}...).

The present paper carries on the computational classification of
closed orientable 3-manifolds performed in \cite{[L]}, by
automatic production and analysis of the complete catalogue of
orientable 3-manifolds represented by crystallizations up to 30
vertices (or, equivalently, admitting coloured triangulations with
at most 30 tetrahedra). It is worth noting that, in the present
work, a suitable use of elementary combinatorial moves yields an
automatic partition of the elements of the generated
crystallization catalogue $\textbf C^{(30)}$ into equivalence
classes, which are proved to be in one-to one correspondence with
the homeomorphism classes of the represented manifolds
(Proposition \ref{biezione}).

If the attention is restricted to prime 3-manifolds not belonging
to the existing Lins's catalogue, the obtained classification may
be summarized by the following statement:

\medskip

\noindent  {\bf Theorem I} \ \ {\it There exist exactly forty-one
closed connected prime orientable 3-manifolds, which admit a
coloured triangulation consisting of 30 tetrahedra and do not
admit a coloured triangulation consisting of less than 30
tetrahedra.

\noindent Among them, there are: \begin{itemize} \item{} 10
elliptic 3-manifolds;
\item{} 17 Seifert non-elliptic 3-manifolds (in particular, 2 torus bundles with Nil geometry);
\item{} 2 torus bundles with Sol geometry;
\item{} 2 manifolds of type $(K \overset{\sim} \times I) \cup (K \overset{\sim}{\times}
I)/A$ \ \ {\rm ($A\in GL(2; \mathbb Z)$, \ $\det(A) = -1$)}, with
Sol geometry;
\item{} 7 non-geometric graph manifolds;
\item{} 3 hyperbolic Dehn-fillings (of the complement of link $6^3_1$).
\end{itemize}}

\bigskip

As a consequence of the generation and analysis of catalogue
$\textbf C^{(30)}$, it is now possible, for any given bipartite
cristallization with at most 30 vertices,\footnote{Non-contracted
graphs representing closed orientable 3-manifolds may be handled
also in case of a higher number of vertices: see Proposition
\ref{riconoscimento}.} to recognize topologically - via computer
program DUKE III\footnote{Details about the C++ program {\it DUKE
III} for automatic analysis and manipulation of PL-manifolds via
edge-coloured graphs may be found on the Web:
http://cdm.unimo.it/home/matematica/ casali.mariarita/DukeIII.htm}
- the represented manifold, with unambiguous identification by
means of JSJ decompositions and fibering structures.

\smallskip

We point out that interesting results follow from a comparative
analysis of both complexity and geometric properties of the
manifolds represented by the subsequent subsets $\mathcal
C^{(2p)},$ \ $1\leq p\leq 15$, of all crystallizations in $\textbf
C^{(30)}$ with exactly $2p$ vertices: in fact, for any fixed
complexity $c$, catalogues $\mathcal C^{(2p)}$ turn out to cover,
for increasing $p,$ first the  most ``complicated" types of
complexity $c$ \ 3-manifolds and then the simplest ones (see Table
2). As a consequence, catalogues $\mathcal C^{(2p)},$ for
increasing value of $p$, appear to be a useful source for
interesting examples in order to test conjectures and search for
patterns about 3-manifolds.

\medskip

Finally, the last paragraph of the paper is devoted to present a
significant improvement of catalogues $\mathcal C^{(2p)}$ (and of
the corresponding catalogues $\tilde {\mathcal C}^{(2p)}$ for
non-orientable 3-manifolds, too): an additional hypothesis on the
representing objects yields a considerable reduction of the
catalogues without loss of generality as far as the represented
3-manifolds are concerned (see Proposition \ref{rigid
cluster-less} and Table 3).

\bigskip
\bigskip
\vskip 2pc
\section{\hskip -0.7cm . Basic notions on coloured triangulations
of manifolds}

\bigskip
As already pointed out, this paper is based on the fundamental
tool of \textit{crystallization theory}, i.e. on the possibility
of representing PL $n$-manifolds by means of \textit{edge-coloured
graphs} or - equivalently - by means of \textit{coloured
triangulations}.

Although crystallization theory extends to manifolds with boundary
and several concepts and results hold for singular manifolds too,
throughout this paper we will restrict our attention to closed,
connected (PL) manifolds.

\bigskip

\noindent \textbf{Definition 1.} An \textit{(n+1)-coloured graph }
is a pair $(\Gamma,\gamma),$ where $\Gamma=(V(\Gamma), E(\Gamma))$
is a regular multigraph\footnote{All notations of general graph
theory are given in accordance to \cite{[W]}.} of degree $n+1$ and
$\gamma : E(\Gamma) \to \Delta_n=\{0,1,\dots,n\}$ is injective on
adjacent edges.

\bigskip

The elements of the set $\Delta_n=\{0,1,\dots,n\}$ are called
\textit{colours}; moreover, for each $i\in\Delta_n$, we denote by
$\Gamma_{\hat i}$ the $n$-coloured graph obtained from
$(\Gamma,\gamma)$ by deleting all edges coloured by
$i$.

An $n$-dimensional pseudocomplex $K$ (see \cite{[HW]} for details)
is called \textit{(vertex)-coloured} if it is equipped with a
labelling of its vertices by $\Delta_n$, which is injective on
each simplex.

The concepts of edge-coloured graph and coloured pseudocomplex are
strictly related; in fact, any $(n+1)$-coloured graph
$(\Gamma,\gamma)$ may be thought of as the dual 1-skeleton of a
coloured $n$-pseudocomplex $K=K(\Gamma)$ (whose $n$-simplices are
in bijection with the vertices of $\Gamma$), so that an
edge-coloration $\gamma$ is naturally induced by that of $K$
(i.e.: for each $e\in E(\Gamma)$, $\gamma(e)=i$ iff the vertices
of the $(n-1)$-simplex of $K$ dual to $e$ are coloured by
$\Delta_n\setminus\{i\}$).

For details about both constructions from edge-coloured graphs to
coloured pseudocomplexes and viceversa, we refer to \cite{[FGG]}
and \cite{[BCG]}.

If polyhedron $|K(\Gamma)|$ is PL-homeomorphic to an $n$-manifold
$M^n$, then $(\Gamma,\gamma)$ is called a \textit{gem} (graph
encoded manifold) of $M^n$, or an edge-coloured graph {\it
representing} $M^n$, while $K=K(\Gamma)$ is said to be a {\it
coloured triangulation} of $M^n$. Furthermore, if
$(\Gamma,\gamma)$ is contracted, i.e. $\Gamma_{\hat i}$ is
connected, for each $i\in\Delta_n$ (equivalently:  $K=K(\Gamma)$
contains only one $i$-coloured vertex, for each $i\in\Delta_n$),
it is called a \textit{crystallization} of $M^n$.

It is easy to see that $M^n$ is orientable iff any edge-coloured
graph $(\Gamma,\gamma)$ representing it is bipartite.

\medskip

Classical results (see \cite{[FGG]}) assure that each $n$-manifold
admits a crystallization; obviously, it generally admits many and
it is a basic problem how to recognize crystallizations (or, more
generally, gems) of the same manifold.

The easiest case is that of two \textit{colour-isomorphic} gems,
i.e. if there exists an isomorphism between the graphs, which
preserves colours up to a permutation of $\Delta_n$. It is quite
trivial to check that two colour-isomorphic gems produce the same
polyhedron.

The following result assures that colour-isomorphic graphs can be
effectively detected by means of a suitably defined numerical
\textit{code}, which can be directly computed on each of them (see
\cite{[BCG]}).

\begin{prop}
Two gems are colour-isomorphic iff  their codes
coincide.
\end{prop}

Also the problem of recognizing non-colour-isomorphic gems
representing the same manifold is solved, but not algorithmically:
a finite set of moves - the so called \textit{dipole moves} - is
proved to exist, with the property that two gems represent the
same manifold iff they can be related by a finite sequence of such
moves (\cite{[FG]}).

However, in this paper, whose results concern dimension three, we
will use another set of moves - \textit{generalized dipole moves}
- defined only for 4-coloured graphs (see section 3).

Even if they still do not solve algorithmically the problem for
general 3-manifolds, nevertheless we will prove that a fixed
sequence of generalized dipole moves is sufficient for classifying
all \textit{rigid} crystallizations of 3-manifolds having at most
30 vertices (see sections 3 and 4).

The definition of \textit{rigid} crystallization requires some
preliminaries.

\bigskip

\noindent \textbf{Definition 2.} A pair $(e,f)$ of distinct
$i$-coloured edges in a 4-coloured graph $(\Gamma,\gamma)$ is said
to form a \textit{$\rho_m$-pair} ($m=2,3$) iff $e$ and $f$ share
exactly $m$ bicoloured cycles of $\Gamma$.

\bigskip

By a {\it $\rho$-pair} we mean a $\rho_m$-pair, for $m\in\{2,3\}$.

\bigskip

\noindent {\bf Definition 3.} A crystallization $(\Gamma,\gamma)$
of a 3-manifold $M^3$ is called \textit{rigid} if it contains no
$\rho$-pairs.

\bigskip

The restriction to the class of rigid crystallizations doesn't
affect the set of represented 3-manifolds, as the following result
proves:

\begin{prop} {\rm (\cite{[C$_1999$]})} \ \label{rigidità} Every closed connected 3-manifold
$M^3$ admits a rigid crystallization. Moreover, if $M^3$ is
handle-free\footnote{$M^3$ contains a handle if it admits a
decomposition $M^3 = J \# H$, where $H$ denotes either the
orientable or non-orientable $\mathbb S^2$-bundle over $\mathbb
S^1$  and $J$ is a suitable non-empty 3-manifold (possibly
homeomorphic to $\mathbb S^3$).} and $(\Gamma,\gamma)$ is any gem
of $M^3$, with $\#V(\Gamma)=2p,$ then there exists a rigid
cristallization $(\bar \Gamma,\bar \gamma)$ of $M^3$, with
$\#V(\bar \Gamma)\le 2p.$
\end{prop}

Hence, in order to obtain a list of all handle-free 3-manifolds
represented by edge-coloured graphs with at most $2p$ vertices
(i.e., according to \cite{[L]} and \cite{[C$_1999$]}, with {\it
gem-complexity} $\le p-1$), it is sufficient to take into account
rigid crystallizations only.

\bigskip

\noindent \textbf{Remark 1 } If $(e,f)$ is a $\rho_3$-pair in a
gem $(\Gamma,\gamma)$ representing a 3-manifold $M^3$ and
$(\Gamma^{\prime},\gamma^{\prime})$ is the gem obtained from
$(\Gamma,\gamma)$ by \textit{switching the $\rho$-pair} (see
Figure 1, or \cite{[L]} for details), then
$(\Gamma^{\prime},\gamma^{\prime})$ represents a 3-manifold $J$ so
that $M^3 = J \# H$, \  $H$ being either the orientable or
non-orientable $\mathbb S^2$-bundle over $\mathbb S^1.$

\bigskip
\smallskip
\centerline{\scalebox{0.6}{\includegraphics{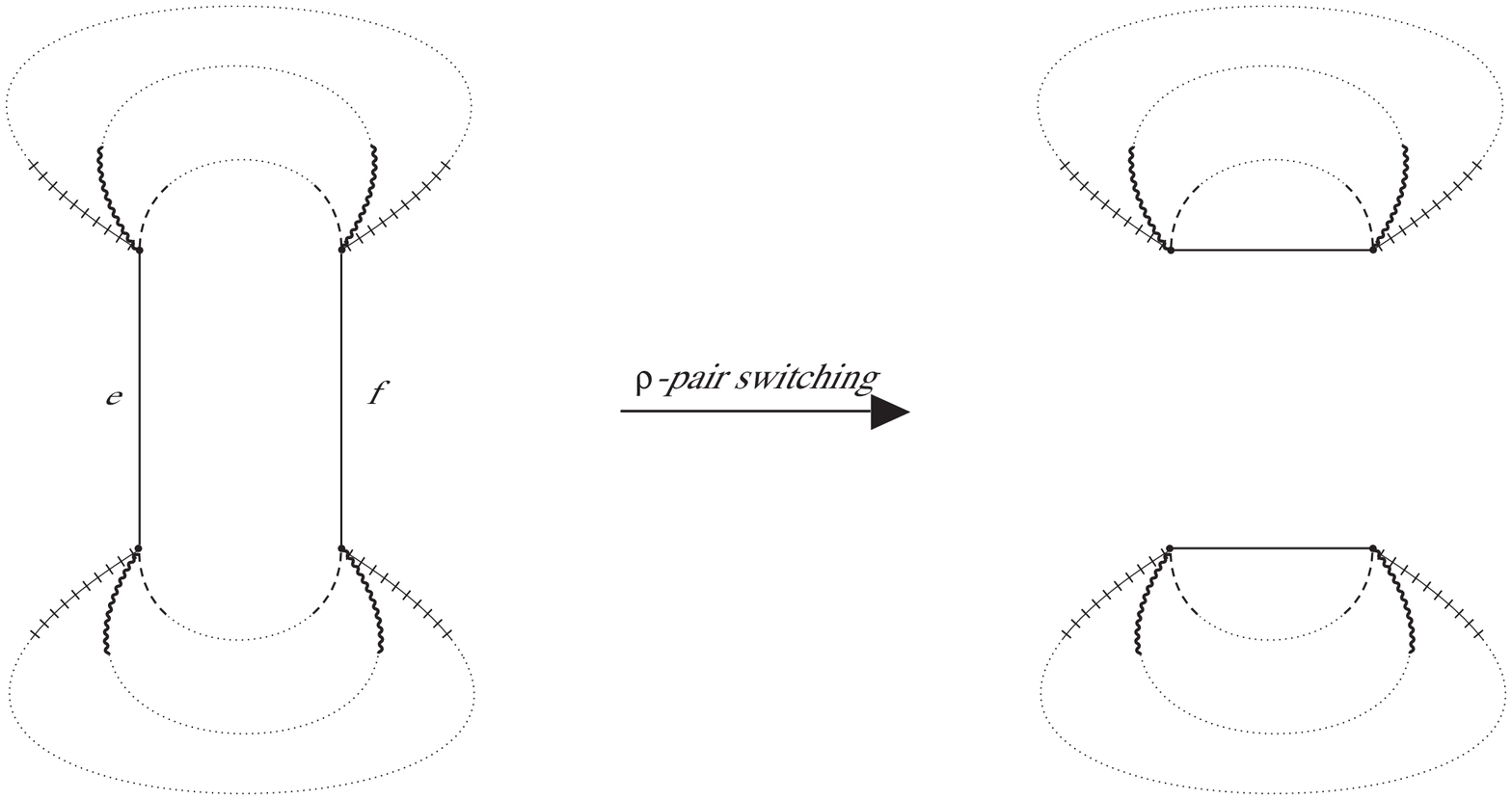}}}
\bigskip
\centerline{{\bf Figure 1}}

\bigskip

Finally, let us point out that an $(n+1)$-coloured graph
$(\Gamma,\gamma)$ is a crystallization of an $n$-manifold  iff
$\Gamma_{\hat i}$ is a gem of $\mathbb S^{n-1}$, for each
$i\in\Delta_n$ (\cite{[FGG]}).

Since a 3-coloured graph represents $\mathbb S^2$ iff it is
planar\footnote{A 3-coloured graph is \textit{planar} iff it has a
cellular embedding in $\mathbb R^2$, whose 2-cells are bounded by
images of bicoloured cycles.}, the characterization of
crystallizations $(\Gamma,\gamma)$ for $n=3$ requires
$\Gamma_{\hat i}$ to be planar and connected for each
$i\in\Delta_n$.

\bigskip

\bigskip

\vskip 2pc

\section{\hskip -0.7cm . Automatic cataloguing and classifying closed 3-manifolds}

\bigskip

Combinatorial encoding of closed 3-manifolds by crystallizations
allows us to construct an essential catalogue of all contracted
triangulations of closed 3-manifolds up to a certain number of
vertices. Moreover, Proposition \ref{rigidità} tells us that we
can restrict our attention to rigid crystallizations; this fact
yields a basic improvement in the direction of the concrete
realization of the catalogue. By using the codes, we can easily
avoid isomorphic graphs, too.

\medskip

For every $p \in\mathbb N$, let $\mathcal C^{(2p)}$ (resp. $\tilde
{\mathcal C}^{(2p)}$) be the catalogue of all non-isomorphic rigid
bipartite (resp. non-bipartite) crystallizations with $2p$
vertices.

The generating algorithm for $\mathcal C^{(2p)}$ and $\tilde
{\mathcal C}^{(2p)}$ was originally described in \cite{[C$_1999$]}
and it consists of the following steps.

\begin{itemize}
\item [Step 1: ] We construct the set $\mathcal S^{(2p)}= \{ \Sigma^{(2p)}_1,
\Sigma^{(2p)}_2, \dots, \Sigma^{(2p)}_{n_p} \}$ of all (connected)
rigid and planar 3-coloured graphs with $2p$ vertices. The
construction makes use of the results of \cite{[L$_1$]} and
\cite{[L]} and is performed by induction on $p$.

\item [Step 2: ] For each $i=1,2,\dots,n_p$, we add to
$\Sigma^{(2p)}_i$ 3-coloured edges in all possible ways to produce
4-coloured graphs, provided that:
\begin{itemize}
\item [-] no vertices belonging to the same bicoloured cycle are
joined (in particular no multiple edges are created \footnote{The
only case where multiple edges are allowed is $p=1$: the order two
3-coloured graph consisting of three multiple edges is rigid and
planar, and obviously gives rise to a rigid order two
crystallization of $\mathbb S^3$  by addition of another
(3-coloured) multiple edge.}) to satisfy the rigidity condition.
\item [-] for each $m\in\{1,\ldots,p\}$, supposing $\,^m\Lambda$
to be a 4-coloured graph (with boundary) obtained from
$\Sigma^{(2p)}_i$  by adding $m < p$ \ 3-coloured edges, the
subgraphs $\,^m\Lambda_{\hat r}$, for every $r\in \{0,1,2\}$, are
planar. This planarity condition can be easily checked since
$$ \,^m\Lambda_{\hat r} \text{\ is \ planar \ \ \ \ iff} \ \ \ \ 2g_{\hat r} -
\,^{\partial}g_{\hat r} = \sum_{i,j \in \Delta_3-\{r\}} \dot
g_{ij} -m$$ where $ 2g_{\hat r}$  (resp. $ \,^{\partial}g_{\hat
r}$) is the number of connected components (resp. of not regular
connected components) of $\,^m\Lambda_{\hat r}$ and $ \dot g_{ij}$
is the number of closed $\{i,j\}$-coloured cycles of
$\,^m\Lambda.$
\item [-] the resulting 4-coloured (regular) graphs $\{
\Gamma^{(2p)}_{i,1}, \Gamma^{(2p)}_{i,2}, \dots,
\Gamma^{(2p)}_{i,m_i}\}$ are crystallizations (i.e.
$\chi(K(\Gamma^{(2p)}_{i,j}))=0$, for each
$j\in\{1,\ldots,m_i\}$).
\end{itemize}

\item[Step 3: ] Let $Y^{(2p)}=
\{\Gamma^{(2p)}_{i,j}\}_{i=1,\ldots,n_p\ \ j=1,\ldots,m_i}$ be the
set of crystallizations arising from Steps 1 and 2. Then, by
computing and comparing the codes and by checking the bipartition
property, we construct the set $X^{(2p)}$ (resp. $\tilde
X^{(2p)}$) consisting of all non-isomorphic bipartite (resp.
non-bipartite) elements of $Y^{(2p)}$.

\item [Step 4: ] The rigidity condition is checked on the elements of $X^{(2p)}$
(resp. $\tilde X^{(2p)}$) and the catalogue ${\mathcal C}^{(2p)}$
(resp. $\tilde {\mathcal C}^{(2p)}$) is obtained. It contains all
rigid bipartite (resp. non-bipartite) crystallizations with $2p$
vertices.
\end{itemize}
\bigskip
The above algorithm was implemented in a C++ program, whose output
data are presented in Table 1 according to the number of vertices.
\bigskip

 \centerline{ \begin{tabular}{|c|c|c|c|}
  \hline    \ & \ & \ & \ \\
  \hfill  {\bf 2p } \hfill &
{\bf  $\# \mathcal S^{(2p)}$ } & {\bf  $\# \mathcal C^{(2p)}$}
& {\bf $\#\tilde {\mathcal C}^{(2p)}$} \\
 \hline    \ & \ & \ & \ \\
  2 & 1 & 1  & 0 \\ 4 &  0 & 0  & 0  \cr 6 & 0 & 0  & 0 \\  8 &
2 & 1  & 0  \\  10 & 0 & 0  & 0   \\ 12 & 1 & 1  & 0 \\  14 & 1 &
1  & 1 \cr 16 & 2 & 3  & 1   \\  18 & 2 & 4  & 1 \\  20 & 8 & 23 &
9 \\ 22 & 8 & 44 & 12   \\  24 & 32 & 262 & 88 \\ 26 & 57 & 1252 &
480 \\ 28 & 185 & 7760 & 2790 \\  30 & 466 & 56912 & 21804
\\
   \ & \ & \ & \       \\
   \hline \end{tabular}}

 \bigskip
\centerline{\bf Table 1:  rigid crystallizations up to 30
vertices.}

\bigskip

Crystallizations of catalogues ${\mathcal C}^{(2p)}$ (resp.
$\tilde {\mathcal C}^{(2p)}$) up to $p=14$ (resp. $p=13$) were
investigated and the related manifolds identified in \cite{[L]}
(resp. \cite{[C$_1999$]}).

In this paper we face the problem of identifying the 3-manifolds
represented by catalogue ${\mathcal C}^{(30)}$ (forthcoming papers
will take into account the same problem for $\tilde {\mathcal
C}^{(28)}$ and $\tilde {\mathcal C}^{(30)}$).

However, our procedure is completely general and, when applied to
${\mathcal C}^{(2p)}$ and $\tilde {\mathcal C}^{(2p)}$ for $p\leq
14$ and $p\leq 13$ respectively, has yielded the known results.

The basic tool is the possibility of subdividing a given set $X$
of crystallizations into subsets (\textit{classes}) such that each
class contains only crystallizations representing the same
manifold.

Obviously, our hope is to obtain classes large enough to coincide
with the topological homeomorphism classes of the manifolds
represented by the elements of $X$.

To fulfill this aim, we need further notions from
crystallization theory.

\medskip

\par \noindent {\bf Definition 4.} Let $(\Gamma,\gamma)$ be a gem of a closed
connected 3-manifold $M^3$.  If there exists an $\{i,j\}$-coloured
cycle of length $m+1$ and a $\{k,l\}$-coloured cycle of length
$n+1$ in $\Gamma$ (with $\{i,j,k,l\}=\Delta_3$) having exactly one
common vertex $\bar v$, then $(\Gamma,\gamma)$ is said to contain
a {\it (m,n)-generalized dipole} of type $\{i,j\}$ at vertex $\bar
v$.
\medskip

To \textit{cancel} a $(m,n)$-generalized dipole from a gem
$(\Gamma,\gamma)$ means to perform on $(\Gamma,\gamma)$ the
operation visualized in Figure 2 (in case $m=3; n=5$).

\bigskip
\smallskip
\centerline{\scalebox{0.55}{\includegraphics{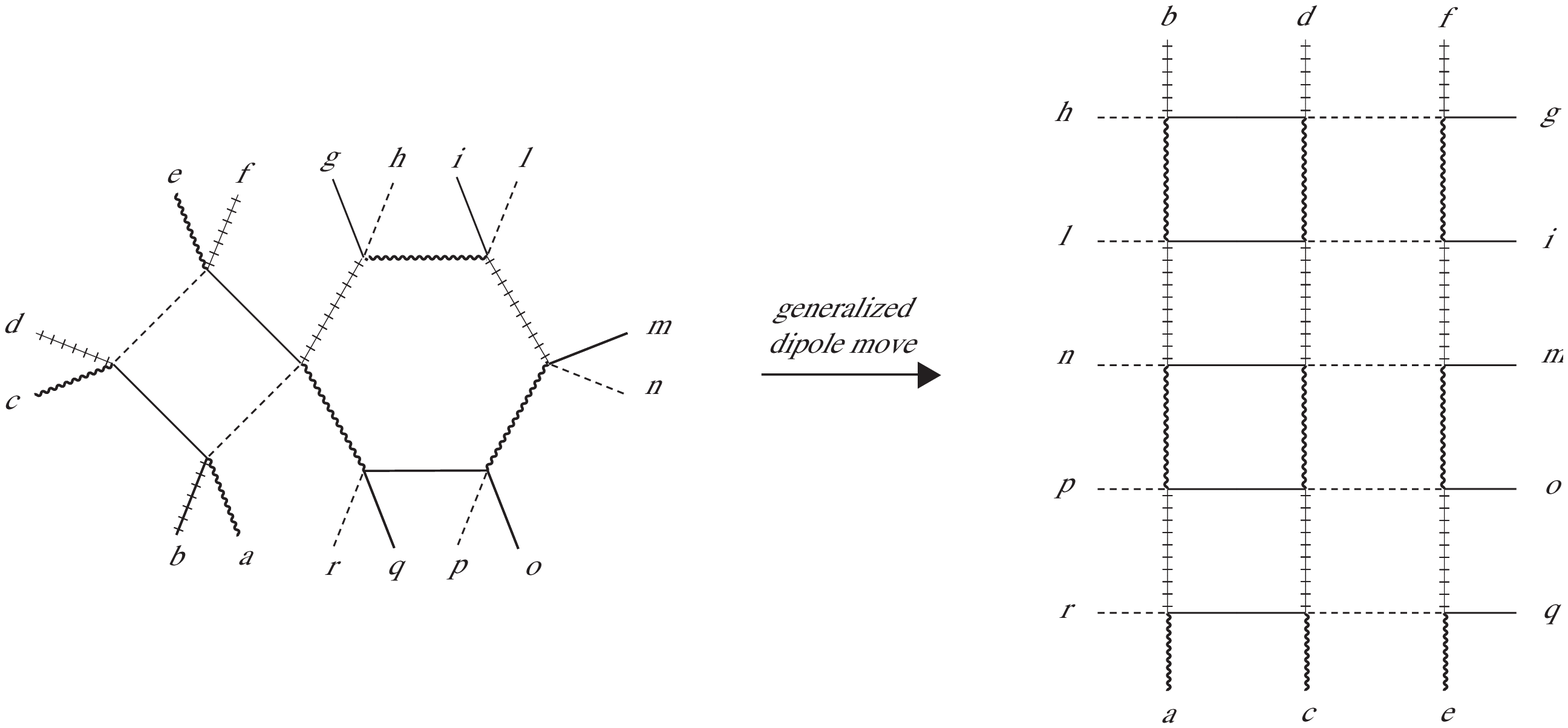}}}
\bigskip
\centerline{{\bf Figure 2}}

\bigskip

In the following we refer to the cancellation of a
$(m,n)$-generalized dipole and to its inverse procedure as
\textit{generalized dipoles moves}.

It is a known result (see \cite{[FG]}) that every two gems which
are transformed into each other by a sequence of generalized
dipoles moves represent the same manifold.

Therefore generalized dipoles moves are a useful tool to
manipulate crystallizations without changing the represented
manifolds.

Let $(\Gamma,\gamma)$ be a rigid crystallization and suppose that in $V(\Gamma)$ an ordering is fixed; given an integer $i\in\{1,2,3\}$, we denote by
$\theta_i(\Gamma)$ the rigid crystallization obtained from
$(\Gamma,\gamma)$ by subsequent cancellations of $(m,n)$-dipoles
of type $\{0,i\},$ according to the following rules:
\begin{itemize}
\item [-] $m,n< 9$ (this condition is necessary
to bound the possible number of vertices of
$\theta_i(\Gamma)).$\footnote{Cancellation of a generalized dipole
increases the number of vertices, but dipoles are frequently
created as a consequence, and their further cancellation allows to
decrease the number of vertices. }

\item [-] supposing $V(\Gamma)=\{v_1,\ldots,v_{2p}\}$, with vertex
labelling coherent with the fixed ordering, the generalized
dipoles of type $\{0,i\}$ are looked for and cancelled for
increasing value of the integer $m\cdot n$ and by starting from
vertex $v_1$ up to $v_{2p}$;  this means that, if $\delta (v_i)$
(resp. $\delta^\prime (v_j)$) is a $(m,n)-$ (resp. a
$(m^\prime,n^\prime)-$) generalized dipole at vertex $v_i$ (resp.
$v_j$), then the cancellation of $\delta (v_i)$ is performed
before the cancellation of $\delta (v_j)$ iff $m\cdot
n<m^\prime\cdot n^\prime,$ or ($m\cdot n=m^\prime\cdot n^\prime$
and $i<j$).

\item [-] after each generalized dipole cancellation,
proper dipoles and $\rho$-pairs, if any, are cancelled in the
resulting graph.
\end{itemize}

Moreover, we define $\theta_0(\Gamma)=(\Gamma,\gamma)$.

\medskip

Note that, given a rigid crystallization $(\Gamma,\gamma)$, there
is an obvious procedure which, starting from the code of
$(\Gamma,\gamma)$, yields a rigid crystallization
$(\Gamma^{<},\gamma^{<})$ which is colour-isomorphic to
$(\Gamma,\gamma)$ and such that an ordering is induced in
$V(\Gamma^{<})$ by the \textit{rooted numbering algorithm}
generating the code (see \cite{[L]}).
As a consequence, for each $i\in\{0,1,2,3\}$, we can
define a map $\theta_i$ on any set $X$ of rigid crystallizations
by setting, for each $(\Gamma,\gamma)\in X$,
$\theta_i(\Gamma)=\theta_i(\Gamma^{<})$, with the ordering of the
vertices induced by the code of $(\Gamma,\gamma)$.

\medskip

Let us define the set
$S_3^0=\{\varepsilon=(\varepsilon_0=0,\varepsilon_1,\varepsilon_2,\varepsilon_3)\
|\ \varepsilon \text{ is a permutation of } \Delta_3 \}$. For each
$\varepsilon\in S_3^0$ and for each $i\in\Delta_3$ we set
$$\theta_{\ll\varepsilon_i\gg}=\theta_{\varepsilon_i}\circ\theta_{\varepsilon_{i-1}}\circ\ldots\circ\theta_{\varepsilon_0}$$
and denote by $h_{\ll\varepsilon_i\gg}(\Gamma)$ the number of
$\rho_3$-pairs which has been deleted when transforming
$\Gamma^{<}$ in $\theta_{\ll\varepsilon_i\gg}(\Gamma^{<})$
(obviously if no $\rho_3$-pair was deleted, we set
$h_{\ll\varepsilon_i\gg}(\Gamma)=0$).

\medskip

Now we are ready to describe the algorithm which, working on a
given list $X$ of rigid crystallizations, produces a partition of
$X$ into equivalence classes, $\{ cl(\Gamma) \ / \ \Gamma \in X
\},$ such that, $\forall \Gamma^{\prime} \in cl(\Gamma)$,
$\Gamma^{\prime}$ and $\Gamma$ represent the same 3-manifold $M$
\textit{up to addition of handles}, i.e. there exist
$h,k\in\mathbb N\cup\{0\}$ such that $|K(\Gamma)|=M\#_h H$ and
$|K(\Gamma^\prime)|=M\#_k H$, where $H=S^1\times S^2$ or
$H=S^1\tilde\times S^2$ according to the bipartition of $\Gamma$
and $\Gamma^\prime$.\footnote{More precisely, $H=S^1\times S^2$
iff $\Gamma$ and $\Gamma^\prime$ are both bipartite or both
non-bipartite.}

The basic idea is that two crystallizations $(\Gamma,\gamma),
(\Gamma^\prime,\gamma^\prime)$ of $X$ belong to the same class iff
there exist $\varepsilon, \mu\in S_3^0$ and $i,j\in\Delta_3$ such
that $\theta_{\ll\varepsilon_i\gg}(\Gamma)$ and
$\theta_{\ll\mu_j\gg}(\Gamma^\prime)$ have the same code.

We consider $X$ as an ordered list and we shall write
$\Gamma\prec\Gamma^\prime$ if $\Gamma$ comes before $\Gamma^\prime$ in
$X$.

For each $(\Gamma,\gamma)\in X$, the construction of $cl(\Gamma)$
is performed via the following algorithm.

\begin{itemize}
\item [Step 1: ] Set $cl(\Gamma)=\{\Gamma\}$ and $h(\Gamma)=0$.
\item [Step 2: ] For each $\varepsilon\in S_3^0$,
$i\in\Delta_3$ and for each $\Gamma^\prime\in X$ with
$\Gamma^\prime\prec\Gamma$, if there exist $\mu\in S_3^0$ and
$j\in\Delta_3$ such that the codes of
$\theta_{\ll\varepsilon_i\gg}(\Gamma)$ and
$\theta_{\ll\mu_j\gg}(\Gamma^\prime)$ coincide, then

\begin{itemize}
\item[$\bullet$] if $h(\Gamma^\prime)-h_{\ll\mu_j\gg}(\Gamma^\prime)\geq
h(\Gamma)-h_{\ll\varepsilon_i\gg}(\Gamma)$, set
$h(\Gamma^{\prime\prime})=k-h(\Gamma)+h_{\ll\varepsilon_i\gg}(\Gamma)+h(\Gamma^\prime)-
h_{\ll\mu_j\gg}(\Gamma^\prime)$ \ \ for each
$\Gamma^{\prime\prime}\in cl(\Gamma)$ with
$h(\Gamma^{\prime\prime})=k;$
\item[$\bullet$] if $h(\Gamma^\prime)-h_{\ll\mu_j\gg}(\Gamma^\prime)<
h(\Gamma)-h_{\ll\varepsilon_i\gg}(\Gamma)$, set
$h(\Gamma^{\prime\prime})=k+h(\Gamma)-h_{\ll\varepsilon_i\gg}(\Gamma)-h(\Gamma^\prime)+
h_{\ll\mu_j\gg}(\Gamma^\prime)$ \ \ for each
$\Gamma^{\prime\prime}\in cl(\Gamma^\prime)$ with
$h(\Gamma^{\prime\prime})=k;$

\end{itemize}
\item [] In both cases, set $c=cl(\Gamma)\cup cl(\Gamma^\prime)$ and
$cl(\Gamma^{\prime\prime})=c$ \ for each $\Gamma^{\prime\prime}\in
c$.
\end{itemize}

\medskip
Furthermore, for each class
$c_i=\{\Gamma_1^i,\ldots,\Gamma_{r_i}^i\}$ and for each $0\leq
h\leq max\{h(\Gamma_1^i), \ldots,$ $h(\Gamma_{r_i}^i)\}$, we
define a partition of $c_i$ into subsets $c_{i,h}=\{\Gamma^i_j\in
c_i\ |\ h(\Gamma^i_j)=h\}$.

\bigskip

Via Proposition \ref{rigidità} and Remark 1 it is very easy to
check that, if $\Gamma\in X$ represents the manifold $M$ with
 $h(\Gamma)=h$ and $c_i=cl(\Gamma)$, then each
 element of $c_{i,k}$ ($0\leq k\leq max\{h(\Gamma^\prime)\ |\ \Gamma^\prime\in c_i\}$) represents
the manifold $M^\prime$ with $M^\prime=M\#_{k-h}H$ or
$M=M^\prime\#_{h-k}H$  (where $H= S^1\times S^2$ or
$H=S^1\tilde\times S^2$, as above), according to $k \ge h$ or
$k<h$.

\bigskip

\noindent\textbf{Remark 2} Note that the algorithm works for any
chosen sequence of generalized dipoles moves; to improve the
results of our implementation, for example, we choose to compare
not only the graphs $\theta_{\ll\varepsilon_i\gg}(\Gamma)$ but
also the graphs
$$\theta_{\ll{\varepsilon^{(k)}_i}\gg}\circ\theta_{\ll{\varepsilon^{(k-1)}_3}\gg}\circ\ldots\circ\theta_{\ll{\varepsilon^{(1)}_3}\gg}(\Gamma),$$
where $\varepsilon^{(k)}$ ($k\in\{2,\ldots,6\}$) is the $k$-th
permutation of $S_3^0$ considered as a lexicographically ordered
set.

\bigskip
\bigskip

In order to analyze and - possibly - identify topologically the
manifolds represented by the crystallizations of a given list $X$,
the first step consists in looking for crystallizations in $X$
which are already identified (as it happens, for example, if known
catalogues of crystallizations are contained in $X$).  In this
case, if the information is added to the input data of our
algorithm, all classes $c_i$ containing at least one of the known
crystallizations, together with their possible subclasses
$c_{i,h}$, turn out to be completely identified: their elements
represent $M\#_t(S^1\times S^2)$ for a fixed manifold $M$ and for
a convenient value of the handle number $t$ (according to the
previously described rules).

\medskip

A second step in the direction of the topological identification
of the given list $X$ is to try connected sums recognition among
manifolds represented by unknown classes. In order to write
explicitly the involved combinatorial condition, we need further
results from crystallization theory.

\begin{prop} \label{somma connessa}
Let $(\Gamma,\gamma)$ be a 3-gem representing a closed connected
3-manifold $M$. Suppose there exist four edges
$\{e_0,e_1,e_2,e_3\}$ in $(\Gamma,\gamma)$ such that
$\gamma(e_i)=i$, for each $i\in\Delta_3$, and
$\Gamma\setminus\{e_0,e_1,e_2,e_3\}$ has two connected components.
Then, there exist 3-gems $(\Gamma^{(1)},\gamma^{(1)})$ and
$(\Gamma^{(2)},\gamma^{(2)})$ such that $M=|K(\Gamma^{(1)})|\#
|K(\Gamma^{(2)})|$.
\end{prop}
\dimo Let $(\bar\Gamma^{(i)},\bar\gamma^{(i)})$, for $i=1,2$, be
the connected components of $\Gamma\setminus\{e_0,e_1,e_2,e_3\}$.
For each $i=1,2$, let us consider the 4-coloured (regular) graph
$(\Gamma^{(i)},\gamma^{(i)})$ obtained from
$(\bar\Gamma^{(i)},\bar\gamma^{(i)})$ by adding a new vertex $v_i$
and four edges $\{e_{i0},e_{i1},e_{i2},e_{i3}\}$ such that
$\gamma^{(i)}(e_{ij})=j$ ($j\in\Delta_3$) and $e_{ij}$ is incident
to $v_i$ and the (boundary) vertex of
$(\bar\Gamma^{(i)},\bar\gamma^{(i)})$ missing colour $j$. It is
easy to see that $(\Gamma^{(i)},\gamma^{(i)})$ is a 3-gem and
$(\Gamma,\gamma)$ is the \textit{connected sum} of
$(\Gamma^{(1)},\gamma^{(1)})$ and $(\Gamma^{(2)},\gamma^{(2)})$
with respect to the vertices $v_1$ and $v_2$ (see \cite{[FGG]}, or
Figure 3 for an example). Therefore $M$ is the connected sum of
$|K(\Gamma^{(1)}|$ and $|K(\Gamma^{(2)})|$ (see \cite{[FGG]}).
\qed

\bigskip
\smallskip
\centerline{\scalebox{0.55}{\includegraphics{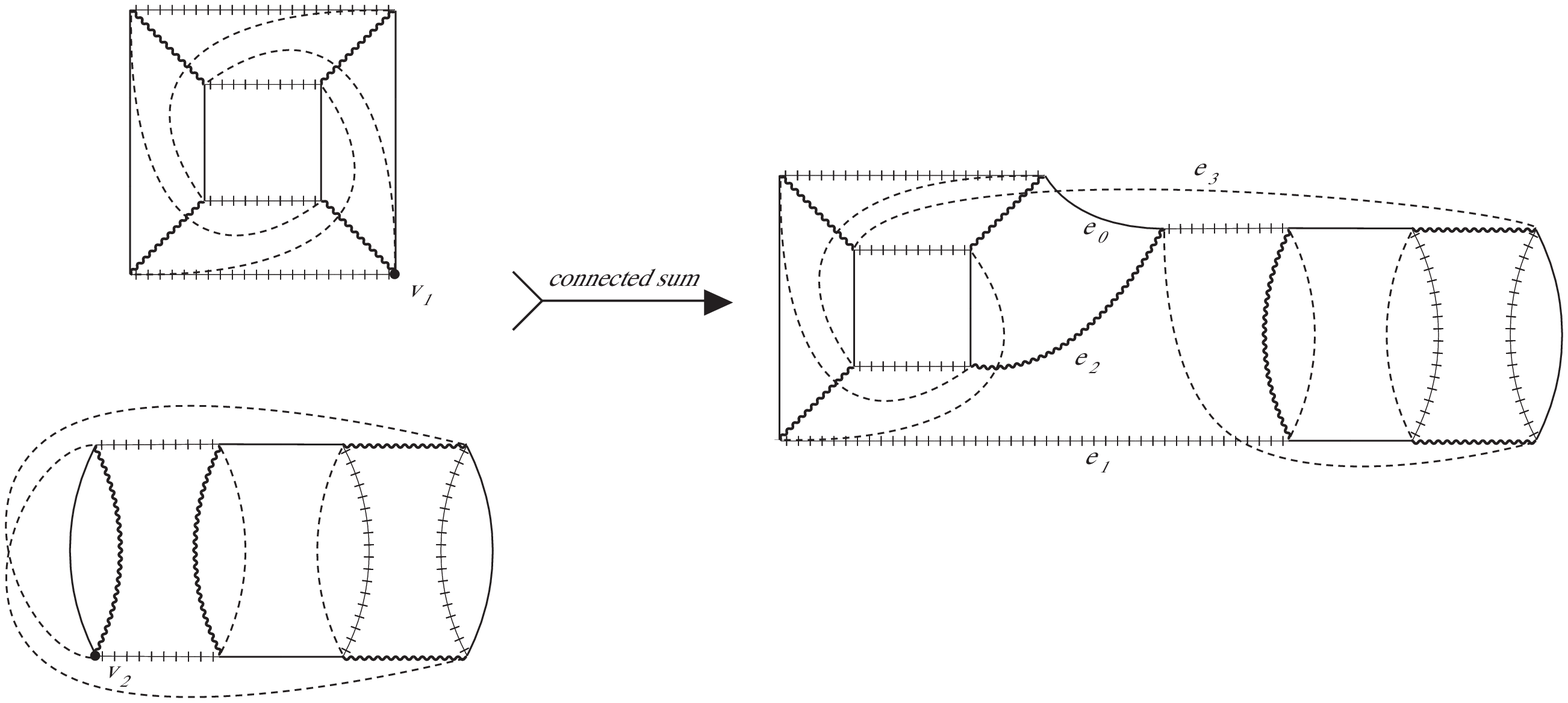}}}
\bigskip
\centerline{{\bf Figure 3}}

\bigskip

Proof of Proposition \ref{somma connessa} tells us that whenever a
crystallization satisfies the condition of the statement we can
split it and analyze the resulting ``pieces".

\bigskip

\noindent\textbf{Remark 3 }Of course, it is possible that one of
the gems $(\Gamma^{(i)},\gamma^{(i)})$ represents $S^3$, i.e. $M$
splits in a trivial connected sum; in any case, the two gems have
fewer vertices than $(\Gamma,\gamma)$ and they will probably be
easier to be recognized.

\bigskip

The C++ program implementing our ``classification" algorithm was
applied to the catalogue $\textbf{C}^{(30)}=\bigcup_{1\leq p\leq
15}\mathcal C^{(2p)}$; it produced 172 classes, 100 of which
completely recognized by means of the existing results about
catalogues $\mathcal C^{(2p)}$, for $1\leq p\leq 14$.

Another C++ program was used to check the condition of Proposition
\ref{somma connessa} on the crystallizations of every class;
whenever the condition was satisfied by a representative
$(\Gamma,\gamma)$ of a class, the crystallizations
$(\Gamma^{(i)},\gamma^{(i)})$ ($i=1,2$) of Proposition \ref{somma
connessa} were constructed and compared, by the code, with the
already known crystallizations of the catalogues $\mathcal
C^{(2p)}$, for $1\leq p\leq 14$.

Sixty-two classes where thus recognized; thirty-one of these connected sums
already appeared in catalogues $\mathcal C^{(2p)}$, for $1\leq
p\leq 14$.

\bigskip

In the following section we will present the analysis of the
above data, which allowed us to obtain the results of Theorem I
and the complete topological classification of the manifolds
encoded by the crystallizations of catalogue $\textbf C^{(30)}$.

\bigskip

\vskip 2pc

\section{\hskip -0.7cm . A complete analysis of catalogue $ \bold{ C^{(30)}}$}

\bigskip

Before discussing our experimental results on $\textbf C^{(30)},$
it is useful to recall the already known topological
identification of the manifolds involved in $\mathcal C^{(2p)}$,
for $1\leq p\leq 14.$

The classification, which is based on the results of \cite{[L]}, was proved in \cite{[CC]} and used for studying a combinatorial concept of complexity
and its relation with Matveev's complexity; the following statement reproduces it in a slightly different form as to suit our present aims.

\begin{prop}  \label{nuova classificazione C^28}

There exist exactly sixty-nine  closed connected prime orientable
3-manifolds, which admit a coloured triangulation consisting of at
most 28 tetrahedra.

Among them, there are:
\begin{itemize}
\item{} $\mathbb S^3$; \item{}  $\mathbb S^2 \! \times \! \mathbb
S^1$; \item{} the six Euclidean orientable 3-manifolds; \item{}
forty-four elliptic 3-manifolds  (in particular, twenty-three lens
spaces); \item{} fifteen Seifert non-elliptic 3-manifolds \par
(more precisely: \begin{itemize} \item[-] four torus bundles with
Nil geometry; \item[-] three manifolds of type $(K \overset{\sim}
\times I) \cup (K \overset{\sim}{\times} I)/A$ \ {\rm (}$A\in
GL(2; \mathbb Z)$, \ $\det(A) = -1${\rm )}, with Nil geometry;
\item[-] seven manifolds with $SL_2(\mathbb R)$ geometry; \item[-]
a further manifold with $Nil$ geometry); \end{itemize} \item{} two
torus bundles with Sol geometry.
\end{itemize}
\end{prop}

\bigskip

The first step towards the identification of the crystallizations
of $\textbf C^{(30)}$ was to compare the results of the
classification program with the known catalogues $\mathcal
C^{(2p)}$ with $1\leq p\leq 14$. This was made by the classifying
program itself.

More precisely, whenever the algorithm produced a crystallization
having less than 30 vertices, the program searched for it in the
known catalogues and gave the resulting name to its class.

In this way 100 classes were recognized; this is exactly the
number of 3-manifolds admitting a rigid crystallization with $2p$
vertices, $1\leq p\leq 14$.

\medskip

Actually, a careful examination of our output data allows to state
the following:

\begin{lemma} \label{vecchie varietà}
The set of classes obtained from $\textbf C^{(30)}$ and containing
at least one crystallization with less than 30 vertices is in
bijective correspondence with the set of 3-manifolds represented
by $\mathcal C^{(2p)},$ \ $1\leq p\leq 14$.
\end{lemma}

\bigskip

Further identifications were obtained by the analysis of the
output of the ``connected sum" program, as described in the above
section.

\smallskip

Let us introduce a notation, which will be useful in the
following.

Let $M,N$ be two closed orientable 3-manifolds; we denote by
$\mathcal C(M,N)$ the set of classes $cl(\Gamma)$ of
crystallizations $(\Gamma,\gamma)\in\textbf C^{(30)}$ such that
$(\Gamma,\gamma)$ satisfies the condition of Proposition \ref{somma
connessa}, with summands $(\Gamma^{(1)},\gamma^{(1)}),\
(\Gamma^{(2)},\gamma^{(2)})$ and
$\{|K(\Gamma^{(1)})|,|K(\Gamma^{(2)})|\}=\{M,N\}$.

\smallskip

We can summarize our results by the following statement.

\begin{lemma}\label{riassunto somme}
\ \
\begin{itemize}
\item [(i)] For each {\rm $(\Gamma,\gamma)\in\textbf C^{(30)}$} satisfying the condition of Proposition \ref{somma connessa} with summands
$(\Gamma^{(1)},\gamma^{(1)})$ and $(\Gamma^{(2)},\gamma^{(2)})$, \
then   $(\Gamma^{(i)},\gamma^{(i)})$ belongs to {\rm
$\textbf{C}^{(28)}=\bigcup_{1\leq p\leq 14}\mathcal C^{(2p)}$},
for each  $i=1,2;$
\item [(ii)] given two closed orientable 3-manifolds $M,N$ such that $\mathcal C(M,N)\neq\emptyset$, we have
\begin{itemize}
\item [-] $\#\mathcal C(M,N)=1$ iff at least one of $M,\ N$ admits orientation-reversing self-homeomorphisms;
\item [-] $\#\mathcal C(M,N)=2$ iff neither $M$ nor $N$ admit orientation-reversing self-homeo- morphisms;
\end{itemize}
\item [(iii)] If $c,c^\prime \in \mathcal C(M,N),$ with $c \ne c^\prime$, then $c$ and $c^\prime$ represent non-homeomorphic manifolds.
\end{itemize}
\end{lemma}

\dimo Statements (i) and (ii) have been deduced directly by the
classification program results. With regard to statement (iii)
more details are required. Let us consider two classes $c$ and
$c^\prime$ as specified above. We have proved that they represent
different manifolds by the following steps:

\begin{itemize}
\item [-] In virtue of statement (ii), the hypothesis of statement (iii) implies $\mathcal C(M,N)=\{c,c^\prime\},$
with both $M$ and $N$ not admitting orientation-reversing
self-homeomor- phisms; moreover, there exists a representative of
$c$ (resp. $c^\prime$), which is a connected sum of a
crystallization of $M$ and a crystallization of $N$;
\item [-] let
$(\Gamma_1,\gamma_1)$ (resp. $(\Gamma_2,\gamma_2)$) be the first
crystallization of catalogue $\textbf C^{(28)}$  representing $M$
(resp. $N$); \item [-] for each $i=1,2$ fix a bipartition on the
set $V(\Gamma_i)$, choose a vertex $v_i\in V(\Gamma_i)$ ($i=1,2$)
and perform the connected sum of graphs $(\Gamma_1,\gamma_1)$ and
$(\Gamma_2,\gamma_2)$ with respect to vertices $v_1$ and $v_2$,
denoting by $\Gamma_1^+\#\Gamma_2^+$ the resulting
crystallization; \item [-] construct the connected sum of
$(\Gamma_1,\gamma_1)$ and $(\Gamma_2,\gamma_2)$ with respect to
$v_1$ and a vertex $w\in V(\Gamma_2)$ such that $v_2$ and $w$
belong to different bipartition classes of $V(\Gamma_2),$ and
denote it by $\Gamma_1^+\#\Gamma_2^-$;
\end{itemize}

It is easy to see that $\Gamma_1^+\#\Gamma_2^+$ and
$\Gamma_1^+\#\Gamma_2^-$ represent non-homeomorphic manifolds (see
\cite{[He]}), which will be denoted by $M^+\# N^+$ and $M^+\#
N^-$.

By applying the classification program to the list formed by
$\Gamma_1^+\#\Gamma_2^+$, $\Gamma_1^+\#\Gamma_2^-$ and the
crystallizations of $c$ and $c^\prime$, we have obtained that
there are exactly two classes $\bar c$ and $\bar c^\prime$ such
that $\bar c$ (resp. $\bar c^\prime$) contains all
crystallizations in $c$ (resp. $c^\prime$)  and one element of
$\{\Gamma_1^+\#\Gamma_2^+, \Gamma_1^+\#\Gamma_2^-\}$ (resp. the
other element of $\{\Gamma_1^+\#\Gamma_2^+,
\Gamma_1^+\#\Gamma_2^-\}$), i.e. $c$ and $c^\prime$ actually
represent the distinct manifolds $M^+\# N^+$ and $M^+\# N^-$. \qed

\noindent \textbf{Remark 4 } With regard to the manifolds $M,N$
such that $\#\mathcal C(M,N)=2$, we point out that in all cases
except one the manifolds involved are $L(3,1)$ and either a lens
space $L(p,q)$ with $(p,q)\in\{(5,1),(7,2),(8,3)\}$ or the
elliptic manifold $S^3/Q_{12}$. The remaining case is the sum of
two copies of $L(4,1)$.

\bigskip

After the comparison with known catalogues $\mathcal C^{(2p)},$
$1\leq p\leq 14,$ and splitting as connected sum, exactly
forty-one unknown classes of crystallizations in $\textbf
C^{(30)}$ turned out to be still unrecognized.

In order to complete the topological identification of all
represented 3-manifolds, a representative for each unknown class
was handled by \textit{Three-manifold Recognizer}\footnote{It is
available on the Web:
http://www.topology.kb.csu.ru/$\sim$recognizer/}, the program
written by V. Tarkaev as an application of the results about
recognition of 3-manifolds obtained by S.Matveev and his research
group.

\bigskip

The output of Matveev-Tarkaev's program proves that all forty-one
classes under examination are topologically distinct and represent
prime manifolds; therefore, by making use also of our former
analysis (see Lemma \ref{vecchie varietà} and Lemma \ref{riassunto
somme}), we can state the following:

\begin{prop} \label{biezione}
\ \
\begin{itemize}
\item [(i)] \textit{There is a bijective
correspondence between the set of equivalence classes obtained by
the classification program and the set of 3-manifolds represented
by {\rm $\textbf C^{(30)}$};}
\item [(ii)] \textit{all connected sums in {\rm $\textbf C^{(30)}$} are identified by the connected sum program, i.e. each class representing a connected sum
contains at least one crystallization satisfying Proposition \ref{somma connessa}.}
\end{itemize}
\end{prop}

\bigskip

With regard to prime 3-manifolds not appearing in $\textbf
C^{(28)}$, the above described analysis of catalogue ${\mathcal
C}^{(30)}$ may be summarized by the following statement, which
directly implies (via Proposition \ref{rigidità}) Theorem I.

\begin{prop} There exist exactly forty-one closed connected prime orientable
3-manifolds, whose minimal coloured triangulation consists of 30
tetrahedra.

\noindent Among them, there are: \begin{itemize} \item{} 10
elliptic 3-manifolds;
\item{} 17 Seifert non-elliptic 3-manifolds (in particular, 2 torus bundles with Nil geometry);
\item{} 2 torus bundles with Sol geometry;
\item{} 2 manifolds of type $(K \overset{\sim} \times I) \cup (K \overset{\sim}{\times}
I)/A$ \ {\rm (}$A\in GL(2; \mathbb Z)$, \ $\det(A) = -1${\rm )},
with Sol geometry;
\item{} 7 non-geometric graph manifolds;
\item{} 3 hyperbolic Dehn-fillings (of the complement of link $6^3_1$).
\end{itemize}
 \end{prop}

\bigskip

\noindent \textbf {Remark 5 } Some of the ``new" prime 3-manifolds
represented by elements of ${\mathcal C}^{(30)}$ can be actually
identified also within crystallization theory, without the aid of
Matveev-Tarkaev's program:
\begin{itemize}
\item[-] the four torus bundles $TB(A)$ (i.e. those obtained with
$A \in  \{ \begin{pmatrix} 1 & 0 \\ 3 & 1 \end{pmatrix},$
$\begin{pmatrix} -1 & 0 \\ 3 &-1 \end{pmatrix},$  $\begin{pmatrix}
-4 & 1 \\ -1 & 0
\end{pmatrix},$ $\begin{pmatrix} 4 & -1 \\ 1 & 0
\end{pmatrix}\} $) have also been recognized by direct construction of the corresponding edge-coloured graphs $\Gamma_{TB}(A)$ (see \cite{[C$_3$]}),
and then by applying, for each $A$, the classification program to
the list formed by  $\Gamma_{TB}(A)$ and the crystallizations of
the only class admitting the group
$\pi_1(TB(A))$ as fundamental group. \item[-] the two manifolds of type $KB(A)=(K
\overset{\sim} \times I) \cup (K \overset{\sim}{\times} I)/A$
(i.e. those obtained with $A \in \left \{ \begin{pmatrix} 1 & -2
\\ -1 & 1 \end{pmatrix}, \begin{pmatrix} -1 & -1 \\ 1 &
2\end{pmatrix}  \right \}$) have also been recognized by direct
construction of the corresponding edge-coloured graphs
$\Gamma_{KB}(A)$ (see \cite{[C$_2005$]}), and then by applying,
for each $A$, the classification program to the list formed by
$\Gamma_{KB}(A)$ and the crystallizations of the only class
admitting the group $\pi_1(KB(A))$ as fundamental group.
\end{itemize}
Moreover:
\begin{itemize}
\item[-]  all 10 elliptic 3-manifolds (i.e.: $S^3/D_{80}$, $S^3/D_{112}$, $S^3/(Q_{28}\times Z_5)$, $S^3/(Q_{32} \times Z_5)$,
$S^3/(P_{48} \times Z_{11})$, $S^3/(P_{48} \times Z_5)$,
$S^3/(P_{48} \times Z_7)$, $S^3/(P_{120} \times Z_{23})$,
$S^3/(P_{120} \times Z_{17})$,  $S^3/(P_{120} \times Z_{13})$) and
the Seifert manifold $SFS(S^2,(2,1),(4,1),(5,2),(1,-1))$, with
$SL_2R$ geometry, have also been recognized directly by means of
homology computation and/or analysis of a presentation of the
fundamental group\footnote{In some cases, GAP program (see
\cite{[Gap]}) has been useful to handle group presentations, by
computation of the corresponding order and/or analysis of
low-index subgroups.}, together with an estimation of the
complexity via GM-complexity (see \cite{[C$_2004$]});
\item[-] two further manifolds (i.e. the Seifert manifolds
$SFS((3,1),(3,2),(3,2),(1,-1))$, with Nil geometry, and
$SFS(K,(2,1)))$, with $SL_2R$ geometry) are easily recognized by
homology computation, together with an estimation of the
complexity via GM-complexity (see \cite{[C$_2004$]}).
\end{itemize}

\bigskip

The complete list of the forty-one closed connected prime
orientable 3-manifolds represented by ${\mathcal C}^{(30)}$ may be
found in Table 2 of \cite{[CC-Web]}. In analogy to the similar
Table 1 of \cite{[CC-Web]}, containing the sixty-nine closed
connected prime orientable 3-manifolds represented by $\textbf
C^{(28)}$ (see Proposition \ref{nuova classificazione C^28} and
\cite{[CC]}), each manifold is identified by means of its JSJ
decomposition and fibering structure, according to Matveev's
description in \cite{[M$_{table11}$]} (see also \cite{[M$_1$]},
\cite{[M$_3$]} and \cite{[M$_4$]}); to make comparison easier, the
position of each manifold within Matveev's table
\cite{[M$_{table11}$]} is also given.

\bigskip

All manifolds involved in $\textbf C^{(30)}$ have also been
detected within Martelli-Petronio censuses of closed irreducible
orientable 3-manifolds up to complexity 10 (see \cite{[MP$_2$]}),
and interesting results followed from a comparative analysis of
both complexity and geometric properties of manifolds represented
by subsequent catalogues $\mathcal C^{(2p)},$ \ $1\leq p\leq 15$.
In fact, the presence of manifolds in $\textbf C^{(30)}$ with
respect to complexity and geometry (see \cite{[MP$_geometry$]} for
cases up to complexity 9, and \cite{[M]} for the last case) may be
summarized in the following table, where the symbol $x/n$ means
that $x$ 3-manifolds appear in $\textbf C^{(30)},$ among the $n$
ones having the appropriate complexity and geometry, and bold
character is used to indicate that \underline{all} manifolds of
the considered type appear in catalogue $\textbf C^{(30)}$:

\bigskip

\hskip -1.75truecm
\begin{tabular}{r|c c c c c c c c c c}{complexity} & 1  & 2 & 3& 4& 5 & 6 & 7 & 8 & 9 & 10 \\
 \hline
 \hline
 lens & {\bf 2/2} &  {\bf 3/3} &  {\bf 6/6} &  {\bf 10/10} &  0/20 & 0/36 & 0/72 & 0/136 & 0/272 & 0/528 \\
 other elliptic & - &  {\bf 1/1} &  {\bf 1/1} &  {\bf 4/4} &  {\bf 11/11} & 14/25 & 0/45 & 0/78 & 0/142 & 0/270 \\
euclidean &  - & - & - & - & - &  {\bf 6/6} & - & - & - & - \\
Nil       &  - & - & - & - & - &  {\bf 7/7} & 3/10 & 0/14 & 0/15 & 0/15 \\
$H^2 \times S^1$  & - & - & - & - & - & - & - & 0/2 & - & 0/8 \\
$SL_2R$  & - & - & - & - & - & - & 13/39 & 0/162 & 0/513 & 0/1416
\\    Sol  & - & - & - & - & - & - & 4/5 & 2/9 & 0/23 & 0/39 \\
non-geometric & - & - & - & - & - & - &  {\bf 4/4} & 1/35 & 2/185 & 0/777 \\
hyperbolic  & - & - & - & - & - & - & - & - &  2/4 & 1/25 \\
 \hline  \hline TOTAL &  {\bf 2/2} &  {\bf 4/4} &  {\bf 7/7} &  {\bf 14/14} &
 11/31 & 27/74 & 24/175 & 3/436 & 4/1154 & 1/3078
\end{tabular}

\medskip
\centerline{\bf Table 2: 3-manifolds involved in $\textbf
C^{(30)}$}

\bigskip
\bigskip

We think worth noting that, for any fixed complexity $c$,
catalogues $\mathcal C^{(2p)}$ cover, for increasing $p,$ first
the  most ``complicated" types of complexity $c$ 3-manifolds and
then the simplest ones: for a detailed analysis on the subject,
see Table 3 of \cite{[CC-Web]}.\footnote{For example note that, as
far as complexity 4 (resp. 5) is concerned, all 10 lens spaces
appear in $\mathcal C^{(28)},$ while all 4 elliptic 3-manifolds
appear in $\mathcal C^{(22)} \cup \mathcal C^{(24)}$ (resp. none
of the 20 lens spaces appear in $\textbf C^{(30)},$ while all 11
elliptic 3-manifolds appear in $\mathcal C^{(24)} \cup \mathcal
C^{(26)}\cup \mathcal C^{(28)}$).}

As a consequence, catalogues $\mathcal C^{(2p)},$ for increasing
value of $p$, appear to be a useful source for interesting
examples to test conjectures and search for patterns about
3-manifolds.

\bigskip
\bigskip

\vskip 2pc

\section{\hskip -0.7cm . Automatic recognition of orientable 3-manifolds}

\bigskip
\bigskip

In this section we want to point out that our approach to the
study of 3-manifolds with low gem-complexity yields not only the
``list" of involved 3-manifolds, but also the ``list" of all
possible coloured triangulations of manifolds with a given number
of tetrahedra.\footnote{A similar point of view may be found in
Burton's works (see \cite{[B$_1$]}, \cite{[B$_2$]},
\cite{[B$_3$]}, \cite{[B$_4$]}) where all minimal triangulations
of manifolds with low complexity are directly constructed.}

In fact, the production and analysis of catalogue ${\mathcal
C}^{(30)},$ with the topological recognition of all represented
3-manifolds (see Theorem I), enables to answer positively the
following general questions, which are ``dual" to each other.

\medskip

\centerline{\it Given a 4-coloured graph with $2p \le 30$
vertices, } \centerline{\it is it possible to say whether it
represents an orientable 3-manifold $M^3$  } \centerline{\it and -
in the affirmative - to recognize $M^3$?}

\medskip

\centerline{\it Given a coloured 3-dimensional pseudocomplex with
 $2p \le 30$ tetrahedra, }
\centerline{\it is it possible to say  whether it represents an
orientable 3-manifold $M^3$ } \centerline{\it and - in the
affirmative - to recognize $M^3$?}

\bigskip

A suitable option of DUKE III program \footnote{See \
http://cdm.unimo.it/home/matematica/casali.mariarita/DukeIII.htm}
answers completely the above questions, for any 4-coloured graph
$(\Gamma, \gamma)$ (resp. any coloured triangulation
$K=K(\Gamma)$): if either the associated matrix $A(\Gamma)$ or the
code $c(\Gamma)$ of $(\Gamma,\gamma)$ is given as input data, then
- in case $|K(\Gamma)|$ is a 3-manifold $M^3$ - the program
identifies $M^3$ within Matveev's catalogue of closed irreducible
orientable 3-manifolds represented by special spines up to
complexity 11 (see \cite{[M$_{table11}$]}), according to Table 1
and Table 2 of \cite{[CC-Web]}.

\begin{prop} \label{riconoscimento} \ Let $(\Gamma, \gamma)$ be any bipartite
4-coloured graph such that  $$\sum_{i,j\in \Delta_3} g_{ij} -
\sum_{i\in \Delta_3} g_{\hat i} =p  \ \ \ \ \ \text{and} \ \ \ \ \
2p - \sum_{i\in \Delta_3} (g_{\hat i} -1) \le 30,$$ where $
g_{\hat i}$  is the number of connected components  of
$\Gamma_{\hat i}$ and $g_{ij}$ is the number of $\{i,j\}$-coloured
cycles of \ $\Gamma.$

\noindent Then, DUKE III program yields, by means of obtained
results about catalogue $\textbf{C}^{(30)},$ the topological
recognition of the closed orientable 3-manifold \
$M^3=|K(\Gamma)|.$
\end{prop}

\dimo First of all, note that condition $\sum_{i,j\in \Delta_3}
g_{ij} - \sum_{i\in \Delta_3} g_{\hat i} =p $  ensures the
polyhedron $|K(\Gamma)|$ to be a PL 3-manifold $M^3.$ Moreover, it
is very easy to check that exactly $\sum_{i\in \Delta_3} (g_{\hat
i} -1)$ subsequent 1-dipole eliminations (each decreasing by two
the number of vertices) may be performed in $(\Gamma, \gamma),$
giving rise to a crystallization of $M^3.$ Hence, the statement
follows directly from the topological identification of the
manifold represented by each element of catalogue
$\textbf{C}^{(30)}$ (see Theorem I and Proposition \ref{nuova
classificazione C^28}). \qed

\bigskip

\vskip 2pc

\section{\hskip -0.7cm . A more essential catalogue for handle-free 3-manifolds}

\bigskip
\bigskip

As already pointed out in paragraph 3, the basic result for the
concrete realization of catalogues $\mathcal C^{(2p)}$ and $\tilde
{\mathcal C}^{(2p)}$ is Proposition \ref{rigidità}, which allows
to restrict the generation process to {\it rigid}
crystallizations, without loss of generality as far as represented
3-manifolds are concerned.

In this paragraph a further, significant improvement in the same
direction is presented; it relies on an idea originally due to
Lins (see \cite{[L]}, paragraph 4.1.4, where basic concepts for
the following definitions and result appear, without successive
application to catalogue generation process yet).

\medskip

\par \noindent {\bf Definition 5.} Let $(\Gamma,\gamma)$ be a gem of the
3-manifold $M^3$. A vertex $v \in V(\Gamma)$ is said to be a {\it
cluster-type vertex} if $(\Gamma,\gamma)$ has four bicoloured
cycles containing $v$ with length four, which involve exactly nine
vertices of $(\Gamma,\gamma).$

\medskip

\par \noindent {\bf Definition 6.} A four-coloured graph $(\Gamma,\gamma)$ is
said to be a {\it cluster-less gem} of a 3-manifold $M^3$ if
$|K(\Gamma)|=M^3$ and $(\Gamma,\gamma)$ admits no cluster-type
vertices.
\bigskip

\begin{prop} \label{clusterless} \ Let $(\Gamma,\gamma)$ be
a gem of a 3-manifold $M^3$, with $V(\Gamma)=2p$. If
$(\Gamma,\gamma)$ contains a cluster-type vertex, then there
exists a cluster-less gem (in particular, a cluster-less
crystallization) $(\Gamma^{\prime},\gamma^{\prime})$ of $M^3$,
with $\#V(\bar \Gamma)< 2p.$ Moreover, if $M^3$ is handle-free and
$(\bar \Gamma,\bar \gamma)$ is a rigid crystallization of $M^3$,
with $\#V(\bar \Gamma)=2\bar p,$ which contains a cluster-type
vertex, then there exists a rigid cluster-less cristallization
$(\bar \Gamma^{\prime},\bar \gamma^{\prime})$ of $M^3$, with
$\#V(\bar \Gamma^{\prime})< 2\bar p.$
\end{prop}

\dimo

As shown in Proposition 24 of \cite{[L]}  (where the hypothesis
that a cluster-type vertex has to involve exactly nine vertices is
actually understood), two cases may arise, for each cluster-type
vertex $v$:
\begin{itemize}
\item{} If the bicoloured cycles with length
greater than four and containing $v$  have a common colour,
 then $(\Gamma,\gamma)$ may be simplified by means of a so called {\it
$TS_1$-move} (which is realized by a standard sequence of dipole
moves, not affecting the number of vertices: see Figure 4, or
Paragraph 4.1.2 of \cite{[L]} for details about $TS$-moves),
followed by a 2-dipole elimination.

\bigskip
\centerline{\scalebox{0.55}{\includegraphics{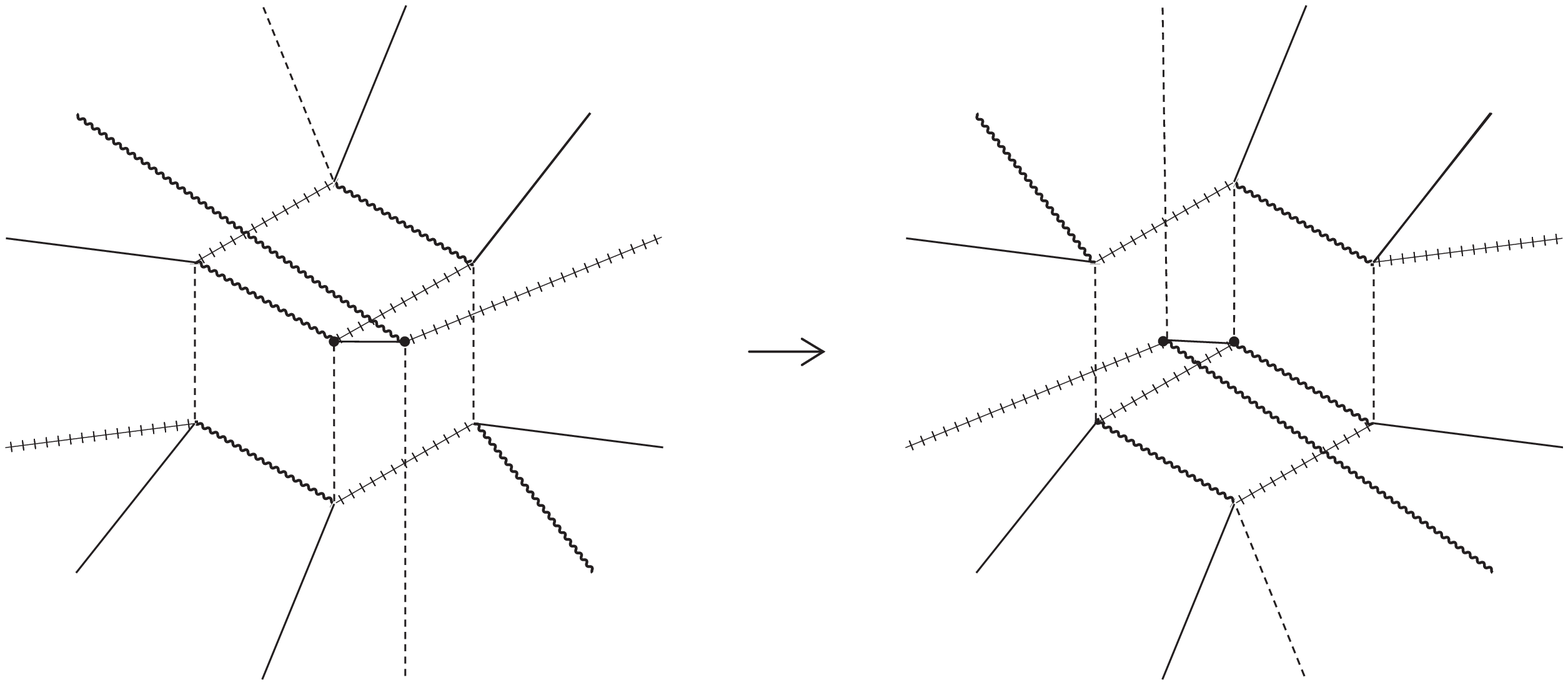}}}
\medskip \centerline{{\bf Figure 4}}

\bigskip

\item{} If the bicoloured cycles with length greater than four and
containing $v$ have no common colours, then $(\Gamma,\gamma)$ may
be thought of as the gem obtained by elimination of a generalized
dipole of type $(3,3)$ on a suitable gem (with exactly two fewer
vertices than $(\Gamma,\gamma)$) of the same manifold (see Figure
5).
\end{itemize}

\centerline{\scalebox{0.75}{\includegraphics{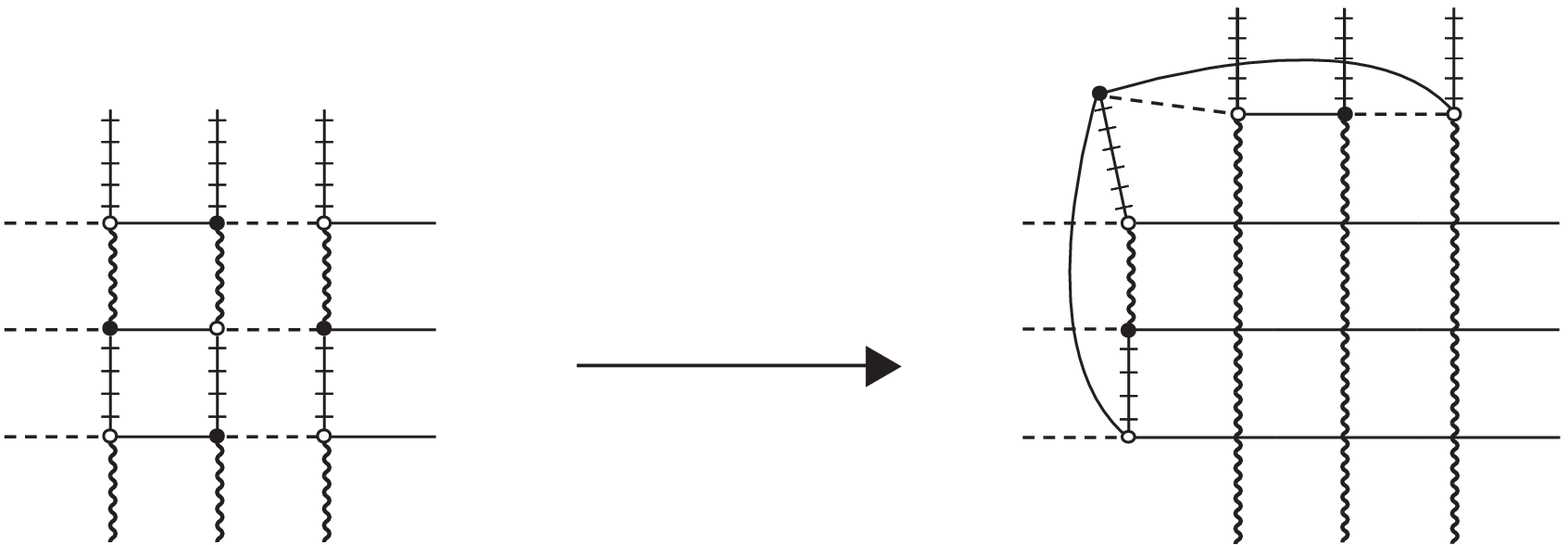}}}
\medskip \centerline{{\bf Figure 5}}

\bigskip

By iterating the process for each cluster-type vertex, a
cluster-less gem $(\Gamma^{\prime},\gamma^{\prime})$ of $M^3$ is
obviously obtained, having strictly fewer vertices than
$(\Gamma,\gamma)$. Moreover, if a cluster-less crystallization of
$M^3$ is required, it is sufficient to perform also all possible
1-dipole eliminations in $(\Gamma^{\prime},\gamma^{\prime}),$ and
then - if it is necessary - to repeat the procedure: since both
eliminations of cluster-type vertices and of 1-dipoles reduce the
number of vertices, the whole process ends in a finite number of
steps, yielding a cluster-less cristallization of $M^3.$

\smallskip

The second part of the statement is a direct consequence of the
second part of Proposition \ref{rigidità}: if $M^3$ is handle-free
and $(\Gamma,\gamma)$ is any gem of $M^3$, both constructions of a
rigid crystallization and of a cluster-less crystallization
subsequently reduce the number of vertices, and hence the
iteration of both processes necessarily yields a rigid
cluster-less cristallization of $M^3.$ \vskip -0.5truecm \ \qed

\smallskip

\begin{prop} \label{rigid cluster-less} \ Every closed connected 3-manifold
$M^3$ admits a rigid cluster-less crystallization.
\end{prop}

\dimo

Let us assume $M^3 = J \#_m H,$ \ $\#_m H$ being the connected sum
of $m \ge 0$ copies of either the orientable or non-orientable
$\mathbb S^2$-bundle over $\mathbb S^1$ and $J$ being a closed
handle-free 3-manifold.  The statement may be easily proved by
making use of the following fundamental facts:
\begin{itemize} \item[-] A rigid cluster-less crystallization of $J$ may be obtained directly by means of Proposition
\ref{clusterless}, applied to any gem of $J$.
\item[-] There exists a rigid cluster-less crystallization
  $(\Omega, \omega)$ (resp. $(\tilde \Omega,\tilde\omega)$) of the
   orientable (resp. non-orientable) $\mathbb S^2$-bundle over $\mathbb
   S^1:$  see Figure 6(a) (resp. Figure 6(b)).
\item[-] If $(\Gamma_1,\gamma_1)$ (resp.
   $(\Gamma_2,\gamma_2)$)  is a rigid cluster-less crystallization of $M^3_1$
  (resp.  $M^3_2$), then the ``connected sum'' $ \Gamma_1 \#
  \Gamma_2$   is a rigid cluster-less
  crystallization of $M^3_1 \# M^3_2$.  In fact, $ \Gamma_1 \#
  \Gamma_2$ is obviously a rigid crystallization (see \cite{[C$_1999$]}, proof of
  Proposition 4), and the absence of cluster-type vertices may be proved easily  by definition of graph connected sum.

\rightline{$\Box$ \ \ \ \ \ \ \ } \end{itemize}

\bigskip
\smallskip
\centerline{\scalebox{0.65}{\includegraphics{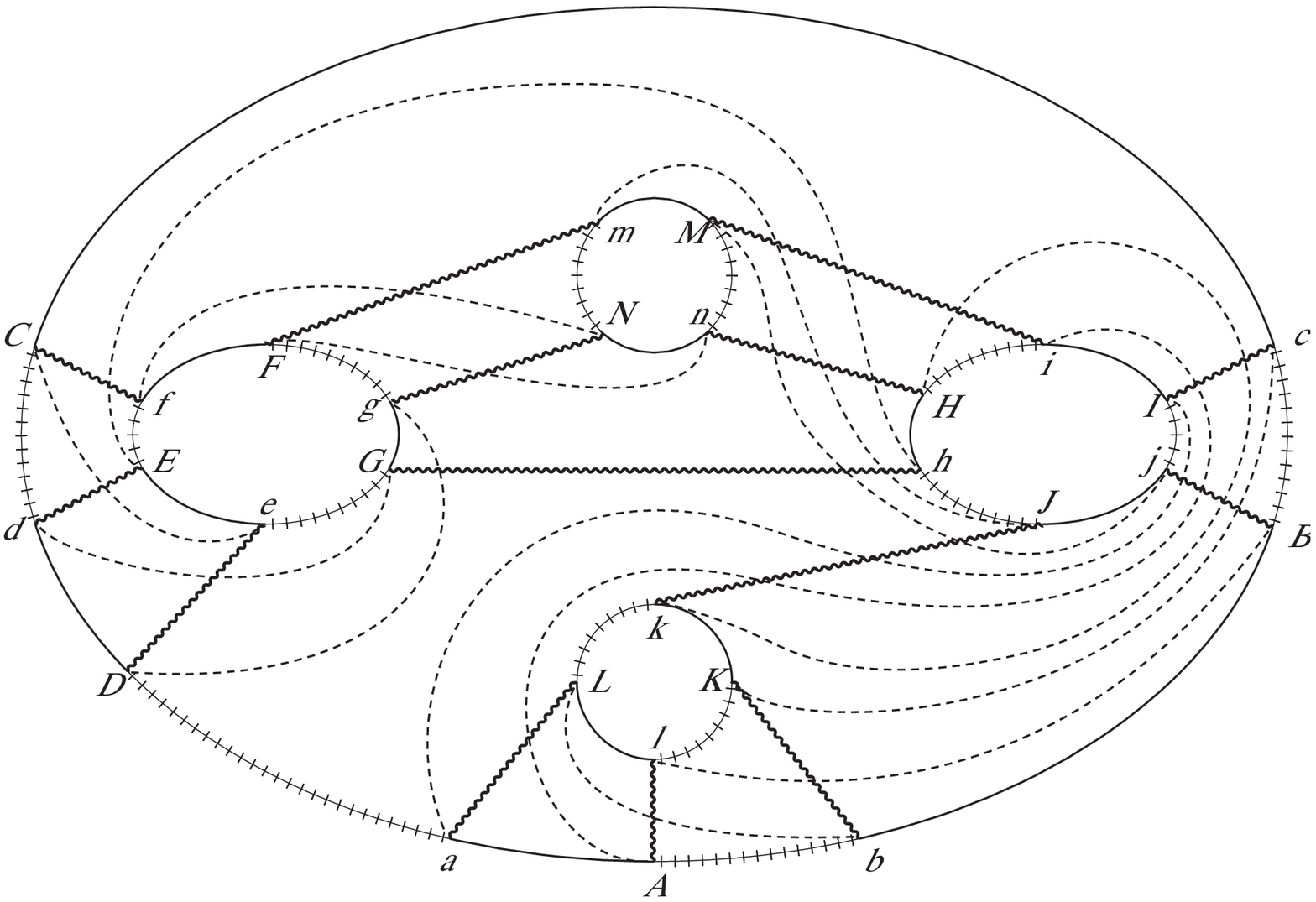}}}
\bigskip
\centerline{{\bf Figure 6(a)}}

\bigskip
\smallskip
\centerline{\scalebox{0.65}{\includegraphics{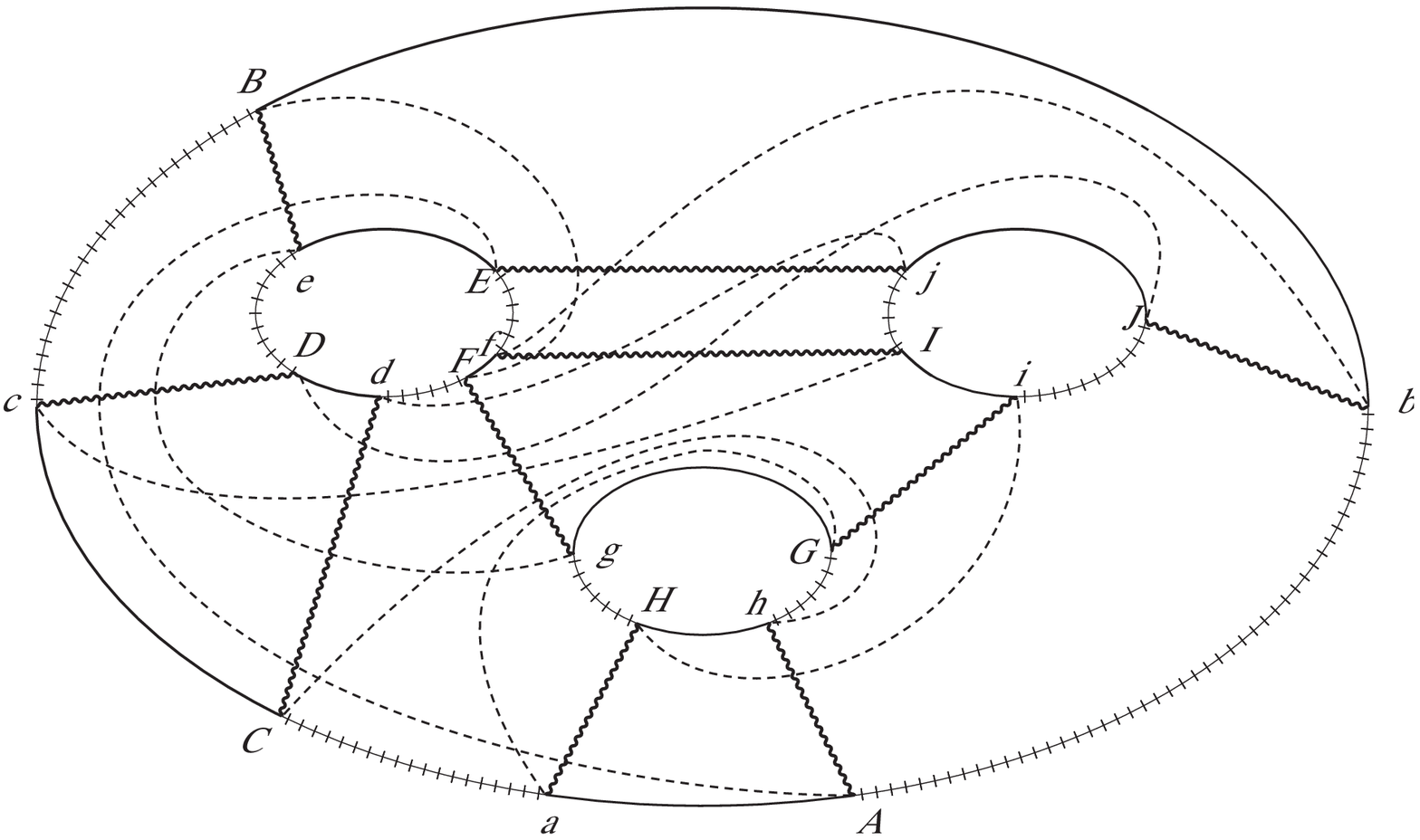}}}
\bigskip
\centerline{{\bf Figure 6(b)}}

\bigskip

\bigskip

Proposition \ref{rigid cluster-less} suggests naturally the
construction of new catalogues $\mathcal C^{\prime \, (2p)}$
(resp. $\tilde {\mathcal C}^{\prime \, (2p)}$) containing rigid
bipartite (resp. non bipartite) cluster-less crystallizations with
$2p$ vertices: in fact, for increasing value of $p$, they yield an
exhaustive representations of {\it all} orientable (resp.
non-orientable) closed 3-manifolds.

\medskip

The following table enables to compare the cardinality of
catalogues  $\mathcal C^{\prime \, (2p)}$ and $\mathcal C^{(2p)}$
(resp. $\tilde {\mathcal C}^{\prime \, (2p)}$ and $\tilde
{\mathcal C}^{(2p)}$), for $1 \le p \le 15$. Notwithstanding the
significative cut in the number of involved crystallizations,
Proposition \ref{clusterless} ensures that, if the attention is
restricted to the handle-free case (as it obviously happens, in
particular, when prime 3-manifolds are considered), no 3-manifold
represented by $\textbf{C}^{(2 \bar p)}=\bigcup_{1\leq p\leq \bar
p}\mathcal C^{(2p)}$  (resp. $\tilde{\textbf{C}}^{(2 \bar
p)}=\bigcup_{1\leq p\leq \bar p}\tilde{\mathcal C}^{(2p)}$), for a
fixed $\bar p$, is lost, when restricting to $\textbf{C}^{\prime
\, (2 \bar p)}=\bigcup_{1\leq p\leq \bar p} \mathcal C^{\prime \,
(2p)}$ (resp. $\tilde{\textbf{C}}^{\prime \, (2 \bar
p)}=\bigcup_{1\leq p\leq \bar p} \tilde{\mathcal C}^{\prime \,
(2p)}$).
\bigskip
\bigskip

\centerline{
  \begin{tabular}{|c|c|c|c|c|}
  \hline    \ & \ & \ & \ & \ \\
  \hfill  {\bf 2p } \hfill &
 {\bf
$\# \mathcal C^{(2p)}$} & {\bf $\# \mathcal C^{\prime \, (2p)}$ }
&
{\bf $\# \tilde {\mathcal C}^{(2p)}$} & {\bf  $\#\tilde {\mathcal C}^{\prime \, (2p)}$} \\
 \hline    \ & \ & \ & \ & \ \\
  2 &  1 & 1  & 0 & 0       \\ 4 &  0 & 0 & 0 & 0 \cr 6 &  0 & 0 & 0 & 0 \\  8  & 1 & 1 & 0 & 0 \\  10 &  0 & 0 & 0  & 0 \\ 12 & 1 & 1 & 0 & 0 \\  14 &
1  & 1 & 1 & 0 \cr 16 & 3 & 3 & 1 & 1 \\  18  & 4 & 2 & 1 & 0 \\
20 & 23 & 16 & 9 & 2 \\ 22 & 44 & 20 & 12 & 4  \\  24 & 262 & 114 & 88 & 17 \\
26 & 1252 & 382 & 480 & 99 \\ 28 & 7760 & 1981 & 2790 & 494 \\  30
& 56912 & 10921 & 21804 & 2989
\\
   \ & \ & \ & \  & \     \\
   \hline \end{tabular}}

 \bigskip
\centerline{\bf Table 3:  rigid and rigid cluster-less
crystallizations up to 30 vertices.}

\bigskip
\medskip

The algorithm described in section 3, applied to catalogue
$\textbf{C}^{\prime \, (30)}=\bigcup_{1\leq p\leq 15} \mathcal
C^{\prime \, (2p)}$, yields exactly 41 classes representing the
prime orientable 3-manifolds described in Theorem I, 69 classes
representing all prime 3-manifolds already contained in
$\textbf{C}^{(28)}=\bigcup_{1\leq p\leq 14}\mathcal C^{(2p)}$ and
63 classes representing non-trivial connected sums.

The bijective correspondence between classes of crystallizations
and manifolds holds in all cases except for connected sums, more
precisely, for manifolds $L(2,1)\#L(2,1)\#$ $L(2,1)$,
$L(2,1)\#L(2,1)\#L(2,1)\#L(2,1)$ and
$L(2,1)\#L(2,1)\#L(2,1)\#L(3,1)$. On the other hand they are
easily recognized through the connected sum program.

Furthermore, manifolds $L(4,1)\#(S^1\times S^2)$,
$L(5,2)\#(S^1\times S^2)$, $S^3/Q_8\#(S^1\times S^2)$ and
$L(2,1)\#L(3,1)\#(S^1\times S^2)$), already appearing in catalogue
$\textbf{C}^{(30)}$, turn out to have no cluster-less
crystallization up to 30 vertices.

However, Proposition \ref{rigid cluster-less} ensures that these manifolds will appear in successive catalogues $\mathcal C^{\prime \,(2p)}$, with $p>15$.

Hence, the generation and analysis of catalogue
$\textbf{C}^{\prime \, (30)}=\bigcup_{1\leq p\leq 15} \mathcal
C^{\prime \, (2p)}$ yield an alternative (and more efficient
\footnote{On the other hand, note that the procedure described in
paragraph 3 turns out to be more useful to identify the manifolds
represented by {\it all} gems up to 30 vertices: see paragraph
4.}) procedure to prove the statement of Theorem I, together with
the existing results collected in Proposition \ref{nuova
classificazione C^28} (see \cite{[L]} and \cite{[CC]}).

\medskip

We hope the catalogues $\mathcal C^{\prime \, (2p)}$ (resp.
$\tilde {\mathcal C}^{\prime \, (2p)}$) can be useful to classify
and recognize topologically  closed orientable (resp.
non-orientable) 3-manifolds admitting coloured triangulations with
$2p$ tetrahedra, with $p>15$.


\bigskip
\bigskip
\bigskip
\bigskip


\begin{thebibliography}{[BCG]}

\bigskip
\bigskip

\bibitem{[AM_1]}  G. Amendola - B. Martelli, {\em Non-orientable 3-manifolds of small
complexity}, Topology And Its Applications, {\bf 133} (2003),
157-178.

\bibitem{[AM_2]}  G. Amendola - B. Martelli, {\em Non-orientable 3-manifolds of
complexity up to 7}, Topology And Its Applications, {\bf 150}
(2005), 179-195.

\bibitem{[BCG]} P.Bandieri - M.R.Casali - C.Gagliardi,
{\em Representing manifolds by crystallization theory:
foundations, improvements and related results}, Atti Sem. Mat.
Fis. Univ. Modena {\bf Suppl. 49} (2001), 283-337.

\bibitem{[BGR]} P.Bandieri - C.Gagliardi - L.Ricci,
{\em Classifying genus two  3-manifolds up to 34 tetrahedra}, Acta
Applicandae Mathematicae {\bf 86}, 267-283.

\bibitem{[BM]} J.Bracho - L.Montejano, {\em The
combinatorics of colored triangulations of manifolds}, Geom.
Dedicata {\bf 22} (1987), 303-328.

\bibitem{[B$_1$]} B.A.Burton, {\em Face pairing graphs and 3-manifold enumeration}, J. Knot Theory Ramifications {\bf 13}(8) (2004),
1057–-1101.

\bibitem{[B$_2$]} B.A.Burton, {\em Structures of small
closed non-orientable 3-manifold triangulations}, to appear in J.
Knot Theory Ramifications, Math.GT/0311113.

\bibitem{[B$_3$]} B.A.Burton, {\em Observations from the 8-tetrahedron non-orientable census}, to appear in Experiment. Math., Math.GT/0509345.

\bibitem{[B$_4$]} B.A.Burton, {\em Enumeration of non-orientable 3-manifolds using face-paring graphs and union-find}, Math.GT/0604584.

\bibitem{[C$_1989$]} M.R.Casali {\em A catalogue of the genus
two 3-manifolds},  Atti Sem. Mat. Fis. Univ. Modena {\bf  37}
(1989), 207-236.

\bibitem{[C$_1999$]} M.R.Casali, {\em
Classification of non-orientable 3-manifolds admitting
decompositions into $\le$ 26 coloured tetrahedra}, Acta Appl.
Math. {\bf 54} (1999), 75-97.

\bibitem{[C$_3$]} M.R.Casali, {\em
Representing and recognizing torus bundles over $\mathbb S^1$},
Boletin de la Sociedad Matematica Mexicana (special issue in honor
of Fico), {\bf 10 (3)}  (2004), 89-106.

\bibitem{[C$_2004$]} M.R.Casali, {\em Computing Matveev's complexity of non-orientable
3-manifolds via crystallization theory}, Topology and its
Applications {\bf 144} (2004), 201-209.

\bibitem{[C$_2005$]} M.R.Casali, {\em
Representing and recognizing 3-manifolds obtained from $I$-bundles
over the Klein bottle}, to appear.

\bibitem{[CC]} M.R.Casali - P. Cristofori, {\em Computing Matveev's complexity via crystallization
theory: the orientable case}, math.GT/0411633.

\bibitem{[CC-Web]} M.R.Casali - P. Cristofori, {\em Archives of closed 3-manifolds with low gem-complexity}, \ available \, from  the \,  Web
\, page \  http://cdm.unimo.it/home/
matematica/casali.mariarita/DukeIII.htm

\bibitem{[Co]} A.Costa {\em Coloured graphs representing
manifolds and universal maps}, Geom. Dedicata {\bf 28} (1988),
349-357.

\bibitem{[FG]}  M.Ferri -  C.Gagliardi, {\em Crystallization
moves}, Pacific J. Math. {\bf 100} (1982), 85-103.

\bibitem{[FGG]} M.Ferri - C.Gagliardi - L.Grasselli, {\em A
graph-theoretical representation of PL-manifolds. A survey on
crystallizations},   Aequationes Math. {\bf 31} (1986), 121-141.

\bibitem{[Gap]} The GAP Group, {\em GAP - Groups, Algorithms and Programming}, Version 4.4, 2005, http://www.gap-system.org

\bibitem{[He]} J.Hempel, {\em 3-manifolds},  Annals of Math.
Studies, {\bf 86}, Princeton Univ. Press, 1976.

\bibitem{[HW]} P.J.Hilton - S.Wylie, {\em An introduction to
algebraic topology - Homology theory},  Cambridge Univ. Press,
1960.

\bibitem{[KL]} L.H.Kauffman - S.Lins {\em Temperley-Lieb
recoupling theory and invariants of 3-manifolds}, Ann. of Math.
Stud. {\bf 134}, Princeton Univ. Press, Princeton (N.J.), 1994.

\bibitem{[L$_1$]} S.Lins, {\em Toward a catalogue of
3-manifold crystallizations}, Atti Sem. Mat. Fis. Univ. Modena
{\bf 33}, (1984), 369-378.

\bibitem{[L]} S.Lins, {\em Gems, computers and attractors for
3-manifolds}, Knots and Everything {\bf 5}, World Scientific,
1995.

\bibitem{[M]} B.Martelli, {\em Complexity of 3-manifolds},  to appear in ``Spaces of Kleinian groups", London Math. Soc. Lecture Notes Ser. 329 (2006),
Math.GT/0405250.

\bibitem{[MP$_1$]} B.Martelli - C.Petronio, {\em Three-manifolds
having complexity at most 9}, Experimental Mathematics {\bf 10
(2)} (2001), 207-236.

\bibitem{[MP$_geometry$]} B.Martelli - C.Petronio, {\em Complexity of geometric three-manifolds}, Geom. Dedicata {\bf 108} (2004), 15-69.

\bibitem{[MP$_2$]} B.Martelli - C.Petronio,
{\em Census 7, Census 8, Census 9, Census 10}, Tables of closed
orientable irreducible 3-manifolds having complexity $c$, $7\le c
\le 10$ , available from the Web page
http://www.dm.unipi.it/pages/petronio/public\_html/files/3D/c9/c9\_census.html

\bibitem{[M$_1$]} S.Matveev, {\em Complexity theory of
three-dimensional manifolds}, Acta Applicandae Math. {\bf 19}
(1990), 101-130.

\bibitem{[M$_3$]} S.Matveev, {\em Recognition and tabulation of three-dimensional manifolds},
Doklady RAS {\bf 400(1)}(2005), 26-28 (Russian;  English trans. in
Doklady Mathematics, {\bf 71} (2005), 20-22).

\bibitem{[M$_4$]} S.Matveev, {\em Tabulation of three-dimensional manifolds},
Uspekhi Mt. Nauk. {\bf 60(4)}(2005), 97-122 (Russian; English
trans. in Russian Math. Surveys {\bf 60(4)}(2005), 673-698).

\bibitem{[M$_{table11}$]} S.Matveev, {\em Table of closed orientable irreducible three-manifolds up to complexity 11},
available from the Web page http://www.topology.kb.csu.ru/
$\sim$recognizer/

\bibitem{[O]} M.A.Ovckinnikov, {\em The table of 3-manifolds of complexity
7}, Preprint Chelyabinsk State University, 1997.

\bibitem{[P]} M.Pezzana, {\em Sulla struttura topologica delle
variet\`a compatte}, Atti Sem. Mat. Fis. Univ. Modena {\bf 23}
(1974), 269-277.

\bibitem{[V]} A.Vince, {\em $n$-graphs}, Discrete Math.
{\bf 72} (1988), 367-380.

\bibitem{[Va]} A.Valverde Colmeiro {\em Grafos coloreados
localmente regulares representando variedades euclideas}, Tesis
Doctoral (UNED - Madrid), 1995.

\bibitem{[W]} A.T.White, {\em Graphs, groups and surfaces}, North
Holland, 1973.

\end{thebibliography}
\end{document}